\documentclass[12pt]{article}

\usepackage[english]{babel}
\usepackage{amsmath}
\usepackage{amssymb,amsfonts,amsthm}
\usepackage{tikz,color,pgf,graphicx,url}
\usetikzlibrary{calc}
\usetikzlibrary{decorations.pathreplacing}
\usepackage{enumerate}     
\usepackage[T1]{fontenc}
\usepackage{xcolor} 
\usepackage{txfonts} 
\usepackage{hyperref}
\hypersetup{
	colorlinks   = true, 
	urlcolor     = black, 
	linkcolor    = black, 
	citecolor   = black 
}

\usepackage[
top=3cm,
bottom=3cm,
inner=2.5cm,      
outer=2.5cm,
footskip=40pt     
]{geometry}

\newtheorem{theorem}{Theorem}[section]
\newtheorem{corollary}[theorem]{Corollary}
\newtheorem{lemma}[theorem]{Lemma}
\newtheorem{proposition}[theorem]{Proposition}
\newtheorem{conjecture}[theorem]{Conjecture}

\newtheorem{claim}{Claim}

\newcommand{\smallqed}{{\tiny ($\Box$)}}
\DeclareMathOperator{\diam}{diam}

\newcommand{\vertex}{\node[vertex]}
\tikzstyle{vertex}=[circle, draw, inner sep=0pt, minimum size=6pt]

\begin{document}

\title{Revisiting $d$-distance (independent) domination in trees and in bipartite graphs}

\author{
Csilla Bujt\'as $^{a,b,}$\thanks{Email: \texttt{csilla.bujtas@fmf.uni-lj.si}}
\and
Vesna Ir\v si\v c Chenoweth $^{a,b,}$\thanks{Email: \texttt{vesna.irsic@fmf.uni-lj.si}}
\and
Sandi Klav\v zar $^{a,b,c,}$\thanks{Email: \texttt{sandi.klavzar@fmf.uni-lj.si}}
\and
Gang Zhang $^{b,d,}$\thanks{Email: \texttt{gzhang@stu.xmu.edu.cn}}
}

\maketitle

\begin{center}
$^a$ Faculty of Mathematics and Physics, University of Ljubljana, Slovenia\\
\medskip

$^b$ Institute of Mathematics, Physics and Mechanics, Ljubljana, Slovenia\\
\medskip

$^c$ Faculty of Natural Sciences and Mathematics, University of Maribor, Slovenia\\
\medskip

$^d$ School of Mathematical Sciences, Xiamen University, China\\
\medskip
\end{center}

\begin{abstract}
The $d$-distance $p$-packing domination number $\gamma_d^p(G)$ of $G$ is the minimum size of a set of vertices of $G$ which is both a $d$-distance dominating set and a $p$-packing. In 1994, Beineke and Henning conjectured that if $d\ge 1$ and $T$ is a tree of order $n \geq d+1$, then $\gamma_d^1(T) \leq \frac{n}{d+1}$. They  supported the conjecture by proving it for $d\in \{1,2,3\}$. In this paper, it is proved that $\gamma_d^1(G) \leq \frac{n}{d+1}$ holds for any bipartite graph $G$ of order $n \geq d+1$, and any $d\ge 1$. Trees $T$ for which $\gamma_d^1(T) = \frac{n}{d+1}$ holds are characterized. It is also proved that if $T$ has $\ell$ leaves, then  $\gamma_d^1(T) \leq \frac{n-\ell}{d}$ (provided that $n-\ell \geq d$), and $\gamma_d^1(T) \leq \frac{n+\ell}{d+2}$ (provided that $n\geq d$). The latter result extends Favaron's theorem from 1992 asserting that $\gamma_1^1(T) \leq \frac{n+\ell}{3}$. In both cases, trees that attain the equality are characterized and relevant conclusions for the $d$-distance domination number of trees derived.
\end{abstract}


\medskip\noindent
{\bf AMS Subj.\ Class.\ (2020)}: 05C69, 05C05

\section{Introduction}
\label{sec:intro}

Let $G = (V(G), E(G))$ be a graph, $S\subseteq V(G)$, let $d$ and $p$ be nonnegative integers, and let  $d(\,\cdot\,, \,\cdot\,)$ denote the standard shortest-path distance. Then $S$ is a {\em $d$-distance dominating set} of $G$ if for every vertex $u\in V(G)\setminus S$ there exists a vertex $w\in S$ such that $d(u,w) \le d$, and $S$ is a {\em $p$-packing} of $G$ if $d(w,w') \ge p+1$ for every two different vertices $w,w'\in S$. The \emph{$d$-distance $p$-packing domination number} $\gamma_d^p(G)$ of $G$ is the minimum size of a set $S$ which is at the same time $d$-distance dominating set and $p$-packing. (If for some parameters $d$ and $p$ such a set does not exist, set $\gamma_d^p(G) = \infty$.)

The $d$-distance $p$-packing domination number was introduced by Beineke and Henning~\cite{Beineke-1994} under the name {\em $(p,d)$-domination number} and with the notation $i_{p,d}(G)$. With the intention of placing it within the trends of contemporary graph domination theory, the notation $\gamma_d^p(G)$ was recently proposed in~\cite{Bujtas-2025+} and we follow it here. 
In~\cite{Bujtas-2025+} it is proved that for every two fixed integers $d$ and $p$ with $2 \le d$ and $0 \le p \leq 2d-1$, the decision problem whether $\gamma_d^p(G) \leq k$ holds is NP-complete for bipartite planar graphs. Several bounds on $\gamma_d^p(T)$, where $T$ is a tree on $n$ vertices with $\ell$ leaves and $s$ support vertices are also proved, including $\gamma_2^0(T) \geq \frac{n-\ell-s+4}{5}$ and $\gamma_d^2(T) \leq \frac{n-2\sqrt{n}+d+1}{d}$, $d \geq 2$. These results improve or extend earlier results from the literature.

In this paper, our focus is on the invariants $\gamma_d^0$ and $\gamma_d^1$. For the first one we will simplify the notation to $\gamma_d$ because it has been investigated under the name of {\em $d$-distance domination number} of $G$ with the notation $\gamma_d(G)$, see the survey~\cite{henning-2020}. We also refer to~\cite{Cz-2022} for algorithmic aspects. For the total version of this concept see~\cite{Davila-2025}. The second invariant $\gamma_d^1$ deals with $d$-distance dominating sets which are $1$-packings. Note that a set of vertices is a $1$-packing if and only if it is an independent set, hence in this case we will say that $\gamma_d^1(G)$ is the {\em $d$-distance independent domination number} of $G$, cf.~\cite{Davila-2017, Gimbel-1996, henning-2020}. 

Meierling and Volkmann~\cite{Meierling-2005}, and independently Raczek, Lema\'nska, and Cyman~\cite{Raczek-2006}, proved that if $d\ge 1$, and $T$ is a tree on $n$ vertices and with $\ell$ leaves, then $\gamma_d(T) \geq \frac{n - d \ell + 2d}{2d+1}$. On the other hand, Meir and Moon~\cite{Meir-1975} proved that if $d \geq 1$ and $T$ is a tree of order $n \geq d+1$, then $\gamma_d(T) \leq \frac{n}{d+1}$. About twenty years later, in 1991, Topp and Volkmann~\cite{Topp-1991} gave a complete characterization of the graphs $G$ with $\gamma_d(G) = \frac{n}{d+1}$. In 1994, Beineke and Henning~\cite{Beineke-1994} proved that if $d \in \{1,2,3\} $ and $T$ is a tree of order $n \geq d+1$, then $\gamma_d^1(T) \leq \frac{n}{d+1}$. Moreover, they closed their paper with the following: 

\begin{conjecture} {\rm \cite{Beineke-1994}}
\label{con:beineke-henning}
If $d\ge 1$ and $T$ is a tree of order $n \geq d+1$, then $\gamma_d^1(T) \leq \frac{n}{d+1}$.
\end{conjecture}

We point out here that in the book's chapter~\cite{henning-2020}, Conjecture~\ref{con:beineke-henning} is stated as~\cite[Theorem 71]{henning-2020} with the explanation that the above-mentioned bound on $\gamma_d(T)$ due to Meir and Moon~\cite{Meir-1975} is proved in such a way, that the $d$-distance dominating set is also independent. Anyhow, in the next section we prove that the bound holds for all bipartite graphs. In Section~\ref{sec:equality cases} we then characterize trees $T$ of order $n$ for which $\gamma_d^1(T) = \frac{n}{d+1}$ holds. In Section~\ref{sec:upper-bounds-leaves}, we prove that if $T$ has $\ell$ leaves, then  $\gamma_d^1(T) \leq \frac{n-\ell}{d}$ (provided that $n-\ell \geq d$), and $\gamma_d^1(T) \leq \frac{n+\ell}{d+2}$ (provided that $n \geq d$). In both cases, the trees that attain the equality are characterized. Using the fact that $\gamma_d(T) \leq \gamma_d^1(T)$, we also derive analogous bounds for $\gamma_d(T)$ and characterize trees attaining those bounds. In particular, if $T$ is a tree with $\ell$ leaves and of order $n\geq d+\ell$, then 
    $$
    \gamma_d(T) \leq \gamma_d^1(T) \leq 
    \begin{cases}
	\frac{n-\ell}{d}, &\text{if~} n<(d+1)\ell ,\\
	\frac{n}{d+1}, &\text{if~} n=(d+1)\ell ,\\
	\frac{n+\ell}{d+2}, &\text{if~} n>(d+1)\ell ,\\
    \end{cases}
    $$
and the upper bounds are best possible. We conclude the paper with a conjecture.  

In the rest of the introduction additional definitions necessary for understanding the rest of the paper are given. For a positive integer $n$ we will use the convention $[n] = \{1,\dots, n\}$. Let $G$ be a graph. The degree of $u\in V(G)$ is denoted by $\deg_G(u)$ or $\deg(u)$ for short. Further, $\diam(G)$ is the diameter of $G$ and $L(G)$ is the set of its leaves, that is, vertices of degree $1$. We call a $d$-distance $p$-packing dominating set of $G$ of size $\gamma_d^p(G)$ a \emph{$\gamma_d^p(G)$-set}. When $G$ is clear from the context, we may shorten it to  \emph{$\gamma_d^p$-set}. A {\it double star} $D_{r,s}$ is a tree with exactly two vertices that are not leaves, with one adjacent to $r \geq 1$ leaves and the other to $s \geq 1$ leaves. When we say that a {\em path $P$ is attached to a vertex} $v$ of a graph $G$, we mean that $P$ is disjoint from $G$ and that we add an edge between $v$ and an end vertex of $P$. 

\section{Bounding $\gamma_d^1$ for bipartite graphs}
\label{sec:conjecture-proved}

For the main result of this section, we first prove the following. 

\begin{theorem}
    \label{thm:bipartite-graphs-gamma_d^1-partition}
    If $d \geq 1$ is an integer and $G$ is a connected bipartite graph of order at least $d+1$, then $V(G)$ can be partitioned into $d+1$ $d$-distance independent dominating sets.
\end{theorem}

\begin{proof}
Set $Z = \diam(G)$. 

If $Z \leq d$, then each vertex is a $d$-distance dominating set of $G$. Since $G$ is bipartite, a required partition of $V(G)$ can be constructed by considering a bipartition $(X,Y)$ of $G$ and partitioning $X$ and $Y$ into $d+1$ parts appropriately. Hence assume in the rest that $Z \geq d+1$.

Let $P$ be a diametrical path of $G$, let $x$ and $y$ be its end-vertices, and root $G$ at $x$. Let $L_i$, $0\le i\le Z$, be the distance levels with respect to $x$, that is, $L_i = \{u\in V(G):\ d(x,u) = i\}$. Consider now the sets
$$S_i = \bigcup_{k\ge 0} L_{k(d+1) + i}, \quad i\in \{0,1,\dots, d\}\,.$$
We claim that $\{S_0, S_1, \dots, S_d\}$ is a partition of $V(G)$ as stated in the theorem.

Since distance levels of a bipartite graph form independent sets and as $d\ge 1$,  each set $S_i$ is independent. Hence it remains to prove that these sets are $d$-distance dominating sets. 

Let $u$ be an arbitrary vertex of $G$ and assume that $u\in L_s$, where $0\le s\le Z$. If $s\ge d$, then there exists a path of length $d$ between $u$ and a vertex from $L_{s-d}$. This already implies that $u$ is $d$-distance dominated by each of the sets $S_i$, $i\in \{0,1,\dots, d\}$. Hence assume in the rest that $s < d$. Then by a parallel argument, $u$ is $d$-distance dominated by each of the sets $S_i$, $i\in \{0,1,\dots, s\}$. It remains to verify that $u$ is $d$-distance dominated by each of the sets $S_i$, $i\in \{s+1,\dots, d\}$. For this sake consider an arbitrary, fixed $t\in \{s+1,\dots, d\}$. Let $Q$ be a shortest $u,y$-path and recall that by our assumption, $d(u,y) \le Z$. Since every edge of $G$ connects two vertices from consecutive distance levels $L_i$, the path $Q$ necessarily contains a vertex $w\in L_t$. We claim that $d(u,w)\le d$. Suppose on the contrary that $d(u,w) > d$. Since $Q$ is a shortest path, $d(w,y) \ge Z-t$. Using these facts together with $t\le d$, we get  
$$ Z \leq d + (Z-t) < d(u,w) + d(w,y) = d(u,y) \le Z\,,$$
which is not possible. We can conclude that $d(u,w)\le d$. This means that  $u$ is $d$-distance dominated by $S_t$ and we are done. 
\end{proof}

In connection with Theorem~\ref{thm:bipartite-graphs-gamma_d^1-partition} we add that Zelinka~\cite{zelinka-1983} proved that if $d \geq 1$ and $G$ is a connected graph of order at least $d+1$, then $V(G)$ can be partitioned into $d+1$ disjoint $d$-distance dominating sets. In this general case, however, the partition need not be into independent sets. 

The following is an immediate consequence of Theorem~\ref{thm:bipartite-graphs-gamma_d^1-partition}.

\begin{corollary}
    \label{cor:bipartite-graphs-upper-bound-gamma_d^1}
    Let $d \geq 1$ be an integer. If $G$ is a bipartite graph of order $n \geq d+1$, then $\gamma_d^1(G) \leq \frac{n}{d+1}$.
\end{corollary}

Corollary~\ref{cor:bipartite-graphs-upper-bound-gamma_d^1} generalizes~\cite[Theorem 71]{henning-2020}. On the other hand, the upper bound in Corollary~\ref{cor:bipartite-graphs-upper-bound-gamma_d^1} may not hold if $G$ is not bipartite. For example, for $n \geq d+2$ and $k \geq 2$, let $G_{n,k,d}$ be the complete graph $K_n$ with $k$ copies of $P_d$ attached to each vertex. Clearly, $|V(G_{n,k,d})| = n (dk+1)$. While $\gamma_d(G_{n,k,d})=n$, a $d$-distance independent domination needs much more vertices, and it is not hard to deduce that $\gamma_d^1(G_{n,k,d}) = 1 + (n-1)k$. As $n \geq d+2$ and $k \geq 2$, we infer that $\gamma_d^1(G_{n,k,d}) > \frac{|V(G_{n,k,d})|}{d+1}$.

\section{Trees that attain equality in Corollary~\ref{cor:bipartite-graphs-upper-bound-gamma_d^1}}
\label{sec:equality cases}

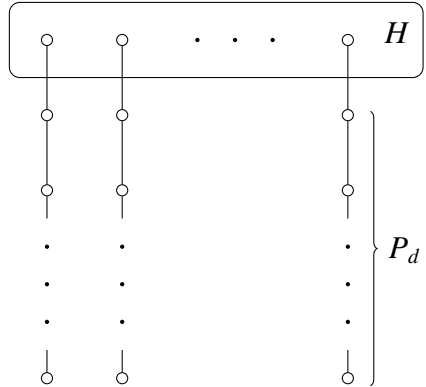
\begin{figure}[ht!]
	\begin{center}
		\begin{tikzpicture}[scale=0.5]
			\tikzstyle{vertex}=[circle, draw, inner sep=0pt, minimum size=6pt]
			\tikzset{vertexStyle/.append style={rectangle}}

  \draw[rounded corners] (0.0, 0.0) rectangle (11.0, 2.0);
                \node ($H$) at (10.3,1.2) {$H$};

                \vertex (1) at (1,1) [scale=0.70,fill=white] {};

                \vertex (1-1) at (1,-1) [scale=0.70,fill=white] {};

                \vertex (1-2) at (1,-3) [scale=0.70,fill=white] {};

                \vertex (1-3) at (1,-4.5) [scale=0.2,fill=black] {};

                \vertex (1-4) at (1,-5.5) [scale=0.2,fill=black] {};

                \vertex (1-5) at (1,-6.5) [scale=0.2,fill=black] {};

                \vertex (1-6) at (1,-8) [scale=0.70,fill=white] {};

                \vertex (2) at (3,1) [scale=0.70,fill=white] {};

                \vertex (2-1) at (3,-1) [scale=0.70,fill=white] {};

                \vertex (2-2) at (3,-3) [scale=0.70,fill=white] {};

                \vertex (2-3) at (3,-4.5) [scale=0.2,fill=black] {};

                \vertex (2-4) at (3,-5.5) [scale=0.2,fill=black] {};

                \vertex (2-5) at (3,-6.5) [scale=0.2,fill=black] {};

                \vertex (2-6) at (3,-8) [scale=0.70,fill=white] {};
                
                \vertex (3) at (9,1) [scale=0.70,fill=white] {};

                \vertex (3-1) at (9,-1) [scale=0.70,fill=white] {};

                \vertex (3-2) at (9,-3) [scale=0.70,fill=white] {};

                \vertex (3-3) at (9,-4.5) [scale=0.2,fill=black] {};

                \vertex (3-4) at (9,-5.5) [scale=0.2,fill=black] {};

                \vertex (3-5) at (9,-6.5) [scale=0.2,fill=black] {};

                \vertex (3-6) at (9,-8) [scale=0.70,fill=white] {};

                \vertex (4) at (5,1) [scale=0.2,fill=black] {};

                \vertex (5) at (6,1) [scale=0.2,fill=black] {};

                \vertex (6) at (7,1) [scale=0.2,fill=black] {};

                \draw[decorate, decoration={brace, raise=5pt}, yshift=-4mm] (9.25,-0.5) -- (9.25,-7.8) node[midway, right=10pt] {};

                \node ($P_d$) at (10.5,-4.6) {$P_d$};

			\path
			(1) edge (1-1)
                (1-1) edge (1-2)
                (1-2) edge (1,-3.75)
                (1,-7.25) edge (1-6)
                (2) edge (2-1)
                (2-1) edge (2-2)
                (2-2) edge (3,-3.75)
                (3,-7.25) edge (2-6)
                (3) edge (3-1)
                (3-1) edge (3-2)
                (3-2) edge (9,-3.75)
                (9,-7.25) edge (3-6)
			;
			
		\end{tikzpicture}
		\caption{The $P_d$-corona $H \circ P_d$ of a graph $H$.}
		\label{fig:path-corona}
	\end{center}
    \end{figure}

Let $d \geq 1$ be an integer. The {\it $P_d$-corona} $H \circ P_d$ of a graph $H$ is the graph obtained from $H$ and $|V(H)|$ disjoint copies of $P_d$, by attaching a copy of $P_d$ to each vertex of $H$, see Fig.~\ref{fig:path-corona}. 

If $d\ge 2$, then let $\mathcal{B}_d$ be the family of $P_d$-coronas of bipartite graphs, that is, 
$$\mathcal{B}_d=\{H \circ P_d: H \text{~ is a bipartite graph}\}.$$
Note that $P_{d+1} \in \mathcal{B}_d$. Observe also that each $G\in \mathcal{B}_d$, where $G = H \circ P_d$, is a bipartite graph with $|V(G)|=(d+1)|V(H)|$. The following proposition shows that the upper bound in Corollary~\ref{cor:bipartite-graphs-upper-bound-gamma_d^1} is best possible.

\begin{proposition}
    \label{prop:extremal-bipartite-graphs-value-gamma_d^1}
    If $G \in \mathcal{B}_d$ is of order $n$, then $\gamma_d^1(G)=\frac{n}{d+1}$.
\end{proposition}

\begin{proof}
    Let $G = H \circ P_d$ for some bipartite graph $H$. By the definition of $H \circ P_d$, the set $L(G)$ is a $d$-distance independent dominating set of $G$. Thus, $\gamma_d^1(G) \leq |L(G)|=|V(H)|=\frac{n}{d+1}$. 
    
    Conversely, for each $u \in V(H)$, let $G_u$ be the subgraph of $G$ induced by $u$ and the vertices of the copy of $P_d$ attached to $u$. Clearly, $G_u \cong P_{d+1}$. If $D$ is a $\gamma_d^1(G)$-set, then $|D \cap V(G_u)| \geq 1$. Thus, $\gamma_d^1(G)=|D| \geq |V(H)|=\frac{n}{d+1}$.
\end{proof}

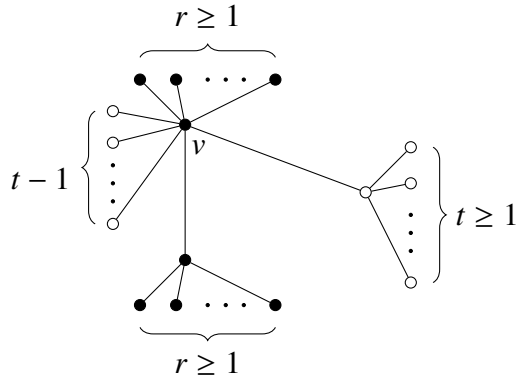
\begin{figure}[ht!]
	\begin{center}
		\begin{tikzpicture}[scale=.6]
			\tikzstyle{vertex}=[circle, draw, inner sep=0pt, minimum size=6pt]
			\tikzset{vertexStyle/.append style={rectangle}}
            
                \vertex (1) at (2,0) [scale=0.7,fill=black] {};
                \vertex (2) at (1,1) [scale=0.7,fill=black] {};
                \vertex (3) at (1.8,1) [scale=0.7,fill=black] {};
                \vertex (4) at (2.5,1) [scale=0.2,fill=black] {};
                \vertex (5) at (2.9,1) [scale=0.2,fill=black] {};
                \vertex (6) at (3.3,1) [scale=0.2,fill=black] {};
                \vertex (7) at (4,1) [scale=0.7,fill=black] {};

                \node (v) at (2.3,-0.5) {$v$};

                \path
                (1) edge (2)
                (1) edge (3)
                (1) edge (7);

                \draw[decorate, decoration={brace, amplitude=5pt, raise=5pt}] (1,1.2) -- (4,1.2) node[midway,above=5pt] {};

                \node (r_geq_1) at (2.5,2.4) {$r \geq 1$};
                
                \vertex (8) at (0.4,0.3) [scale=0.7,fill=white] {};
                \vertex (9) at (0.4,-0.4) [scale=0.7,fill=white] {};
                \vertex (10) at (0.4,-0.9) [scale=0.2,fill=black] {};
                \vertex (11) at (0.4,-1.3) [scale=0.2,fill=black] {};
                \vertex (12) at (0.4,-1.7) [scale=0.2,fill=black] {};
                \vertex (13) at (0.4,-2.2) [scale=0.7,fill=white] {};

                \path
                (1) edge (8)
                (1) edge (9)
                (1) edge (13);

                \draw[decorate, decoration={brace,mirror,amplitude=5pt, raise=5pt}] (0.3,0.3) -- (0.3,-2.2) node[midway,left=5pt, rotate=90] {};

                \node (t_minus_1) at (-1.2,-1.2) {$t-1$};

                \vertex (14) at (2,-3) [scale=0.7,fill=black] {};
                \vertex (15) at (1,-4) [scale=0.7,fill=black] {};
                \vertex (16) at (1.8,-4) [scale=0.7,fill=black] {};
                \vertex (17) at (2.5,-4) [scale=0.2,fill=black] {};
                \vertex (18) at (2.9,-4) [scale=0.2,fill=black] {};
                \vertex (19) at (3.3,-4) [scale=0.2,fill=black] {};
                \vertex (20) at (4,-4) [scale=0.7,fill=black] {};

                \path
                (1) edge (14)
                (14) edge (15)
                (14) edge (16)
                (14) edge (20);

                \draw[decorate, decoration={brace,mirror, amplitude=5pt, raise=5pt}] (1,-4.2) -- (4,-4.2) node[midway,above=5pt] {};

                \node (r_geq_1) at (2.5,-5.3) {$r\geq 1$};

                \vertex (21) at (6,-1.5) [scale=0.7,fill=white] {};
                \vertex (22) at (7,-0.5) [scale=0.7,fill=white] {};
                \vertex (23) at (7,-1.3) [scale=0.7,fill=white] {};
                \vertex (24) at (7,-2) [scale=0.2,fill=black] {};
                \vertex (25) at (7,-2.4) [scale=0.2,fill=black] {};
                \vertex (26) at (7,-2.8) [scale=0.2,fill=black] {};
                \vertex (27) at (7,-3.5) [scale=0.7,fill=white] {};

                \path
                (1) edge (21)
                (21) edge (22)
                (21) edge (23)
                (21) edge (27);

                \draw[decorate, decoration={brace,amplitude=5pt, raise=5pt}] (7.2,-0.5) -- (7.2,-3.5) node[midway,right=5pt, rotate=90] {};

                \node (t_geq_1) at (8.7,-2) {$t \geq 1$};

	    \end{tikzpicture}
		\caption{A tree $T$ from the family $\zeta_1$, where $T'$ is the double star induced by the black vertices.}
		\label{fig:zeta-graph}
	\end{center}
\end{figure}

Note that if $G$ is a connected bipartite graph of order $n=d+1$, then $\gamma_d^1(G)=1=\frac{n}{d+1}$, and $\gamma_d^1(C_{2d+2})=2=\frac{2d+2}{d+1}$. Moreover, if $d=1$, then $\gamma_1^1(K_{r,r})=r=\frac{r+r}{2}=\frac{n}{2}$. In 2004, Ma and Chen gave an equivalent description of the bipartite graphs $G$ of order $n$ with $\gamma_1^1(G)=\frac{n}{2}$, see~\cite[Theorem 1]{MaDeXing-2004}. They also proved an explicit characterization of such a family for the case of trees. To state the result, let $\zeta_1$ be a family of trees defined by the following recursive construction.
\begin{enumerate}
\item[(i)] $K_2\in \zeta_1$.
\item[(ii)] If $T' \in \zeta_1$, and $T$ is obtained by joining the center of a new copy of $K_{1,t}\ (t\geq 1)$ to a support vertex $v$ of $T'$ and adding $t-1$ leaves at $v$, then $T \in \zeta_1$. 
\end{enumerate}
See Fig.~\ref{fig:zeta-graph} for an illustration.\footnote{For the definition of $\zeta_1$, we note that both vertices of a $K_2$ are support vertices. Then, (ii) can be applied to $K_2$, and this step results in the double star $D_{r,r}$ for every $r\ge 1$. It shows that the family $\zeta_1$ is the same as $\{K_2\} \cup \zeta$ in~\cite{MaDeXing-2004}.}


The result of Ma and Chen for trees now reads as follows. 

\begin{theorem}{\rm (\cite[Corollary 1]{MaDeXing-2004})}
    \label{thm:extremal-trees-gamma_1^1}
    If $T$ is a tree of order $n$, then $\gamma_1^1(T)=\frac{n}{2}$ if and only if $T \in  \zeta_1$.
\end{theorem}

We shall focus on the general case for $d \geq 2$, and give a complete characterization of the trees achieving equality in the upper bound of Corollary~\ref{cor:bipartite-graphs-upper-bound-gamma_d^1}. Set
\begin{equation*}
	\mathcal{T}_d=\{T^* \circ P_d: T^* \text{~ is a non-trivial tree}\}.
\end{equation*}
Note that $\mathcal{T}_d$ does not contain the path $P_{d+1}$. Since $\mathcal{T}_d \subset \mathcal{B}_d$, and by Proposition \ref{prop:extremal-bipartite-graphs-value-gamma_d^1}, $\gamma_d^1(T)=\frac{n}{d+1}$ for each tree $T \in \mathcal{T}_d$ of order $n$. Moreover, if $T$ is a tree of order $n=d+1$, then we also have $\gamma_d^1(T)=1=\frac{n}{d+1}$.

In a tree $T$ and for a vertex $v \in V(T)$, let $L(v)$ be the set of leaves of $T$ that are neighbors of $v$ in $T$.  Root $T$ at some vertex. Let $T_v$ be the subtree induced in $T$ by $v$ and its descendants, and let $T-T_v=T-V(T_v)$. A vertex of $T$ is called a \emph{$P_d$-support vertex} if it is attached to a copy of $P_d$. For each $H \in \mathcal{T}_d$, every vertex of $H^*$ is a $P_d$-support vertex of $H$, where $H =H^* \circ P_d$ for some non-trivial tree $H^*$. In particular, a $P_1$-support vertex of $T$ is just a support vertex of $T$. A vertex of $T$ is a \emph{$(P_i, P_j)$-support vertex} if both a copy of $P_i$ and a copy of $P_j$ are attached to it. In particular, a $(P_i,P_i)$-support vertex has at least two copies of $P_i$ attached. The \emph{$d$-subdivision} of $T$ is the tree obtained from $T$ by subdividing each edge $d$-times. Then the $1$-subdivision of $T$ is just the subdivision of $T$.

Before proving the announced characterization of trees of order $n$ with $\gamma_d^1(T)=\frac{n}{d+1}$, we state the following lemma which will also be used in the subsequent section. 

\begin{lemma}
        \label{lemma:structural}
       Let $d \ge 2$ and let $T$ be a tree with $s=\diam(T) \geq 2d+1$. Suppose that $P:=v_1v_2 \dots v_{s+1}$ is a diametrical path in $T$ and the tree is rooted at $v_{s+1}$. If there is no $P_{d+1}$-support vertex and no $(P_i,P_j)$-support vertex in $T$ with $i \in [d-1]$ and $j \in [d]$, then the following statements hold.      
        \begin{itemize}        
            \item[(i)] If $k \in \{2, \dots, d\} \cup \{s-d+2, \dots ,s\}$, then $\deg(v_k) = 2$. \item[(ii)] If $k \in \{ d+1, s - d + 1 \}$, then $\deg(v_{k}) \geq 3$.
            \item[(iii)] For every $v \in V(T)$, if $v$ is the only vertex with $\deg_T(v) \ge 3$ in the subtree $T_v$, then $T_v$ is isomorphic to the $(d-1)$-subdivision of a star $K_{1,t}$ with $t \geq 2$. 
             \item[(iv)] The subtree $T_{v_{d+1}}$ is isomorphic to the $(d-1)$-subdivision of a star $K_{1,t}$ with $t \geq 2$. 
            \item[(v)]  If $s=2d+1$, then $T$ is obtained by taking the $(d-1)$-subdivisions of two stars  $K_{1,t_1}$ with $t_1 \geq 2$ and $K_{1,t_2}$ with $t_2 \geq 2$, and adding an edge between the centers.  
        \end{itemize}           
       \end{lemma}
       \proof
(i)--(ii) According to the conditions,  there is no $(P_1, P_1)$-support vertex in $T$.  That is, every vertex of $T$ is adjacent to at most one leaf, and in particular, $\deg(v_2) = 2$. Further, if $d \ge 3$, then $\deg(v_3) =2$, since otherwise $v_3$ would be a $(P_i,P_2)$-support vertex with $1 \leq i \leq 2 \leq d-1$ contradicting the condition. Similarly, $\deg(v_k) = 2$ holds for all $k \in \{2, \dots, d\}$. By symmetry, the same is true for $v_k$ if $k \in \ \{s-d+2, \dots ,s\}$. This proves (i). The assumption that there is no pendant $P_{d+1}$ in $T$ directly implies (ii).


(iii) If $\deg_T(v) \ge 3$ and $\deg_{T_v}(u) \leq 2$ for every further vertex $u$ from $V(T_v)$, then at least two pendant paths are attached to $v$. Then, by the conditions in the lemma, every path attached is isomorphic to $P_d$.

(iv)--(v) As $P$ is a diametrical path, a vertex $u \in V(T_{v_{d+1}})$ different from $v_{d+1}$ cannot be a $P_d$-support vertex. Part (iii) then implies (iv). If we re-root $T$ at the vertex $v_1$, the same property holds for the subtree induced by $v_{s-d+1}$ and its descendants in the re-rooted tree. This directly implies (v) for the case of $s=2d+1$. \qed
\medskip

\begin{theorem}
    \label{thm:extremal-trees-gamma_d^1}
    If $d \geq 2$ and $T$ is a tree of order $n$, then $\gamma_d^1(T)=\frac{n}{d+1}$ holds if and only if $n=d+1$ or $T \in \mathcal{T}_d$.
\end{theorem}

\begin{proof}
     If $T$ is of order $n=d+1$, then, clearly, $\gamma_d^1(T)=1=\frac{n}{d+1}$, and if $T \in \mathcal{T}_d$, then by Proposition~\ref{prop:extremal-bipartite-graphs-value-gamma_d^1}, we have $\gamma_d^1(T)=\frac{n}{d+1}$. The proof of the necessity is by induction on $n$. If $\gamma_d^1(T)=\frac{n}{d+1}$, then $n=(d+1)q$ for some integer $q \geq 1$. If $q=1$, then $n=d+1$. So, we may assume that $q \geq 2$ and $n \geq 2(d+1)$. If $\diam(T) \leq 2d$, then $\gamma_d^1(T)=1 < \frac{2(d+1)}{d+1} \leq \frac{n}{d+1}$. In the continuation, we assume that $\diam(T) \geq 2d+1$ and $\gamma_d^1(T)=\frac{n}{d+1}$. 
    
    \begin{claim}
        \label{claim:path-support-vertex}
        If $i \in [d-1]$ and $j \in [d]$, then there is no $(P_i, P_j)$-support vertex in $T$.  
    \end{claim}

\proof
Suppose, to the contrary, that $v$ is a $(P_i, P_j)$-support vertex in $T$ and $i \leq j$. Let $P':=x_1x_2\dots x_i$ and $P'':=y_1y_2\dots y_j$ be two copies of $P_i$ and $P_j$ attached to $v$ in $T$, where $x_iv,y_jv \in E(T)$. Note that $d(x_1,v)=i \leq  j = d(y_1,v)$. Consider $T'=T-V(P')$. Then $n'=|V(T')|=n-i  \geq 2(d+1)-(d-1)=d+3$. Let $D'$ be a $\gamma_d^1(T')$-set. If $v \in D'$, then $D'$ is also a $d$-distance independent dominating set of $T$. If $v \notin D'$, then  $|D'|=\gamma_d^1(T')$ implies $|D' \cap V(P'')| \leq 1$. For the subcase $|D' \cap V(P'')| = 1$, we may assume that $y_j \in D'$. Then $d(x_k,y_j) \leq d$ for each $k \in [i]$. For the subcase $|D' \cap V(P'')| = 0$, in order to $d$-distance dominate $y_1$ in $T'$, there exists a vertex $u \in D'$ such that $d_{T'}(u,y_1) \leq d$. Since $i \le j$, it holds that $d_T(x_k,u)\leq d_T(u,y_1)=d_{T'}(u,y_1) \leq d$ for each $k \in [i]$. Thus, $D'$ is always a $d$-distance  independent dominating set of $T$. Corollary~\ref{cor:bipartite-graphs-upper-bound-gamma_d^1} then implies $$\gamma_d^1(T) \leq |D'|=\gamma_d^1(T') \leq \frac{n'}{d+1}
<\frac{n}{d+1},$$
        which contradicts the assumption  $\gamma_d^1(T) = \frac{n}{d+1}$. 
        This proves Claim~\ref{claim:path-support-vertex}.
\smallqed

\medskip
    \begin{claim}
        \label{claim:P_d-support-vertex}
        If $T$ has a $P_{d+1}$-support vertex $v$, then $T \in \mathcal{T}_d$.
    \end{claim}

\proof
Let $P':=x_1x_2\dots x_{d+1}$ be a copy of $P_{d+1}$ attached to $v$, where $x_{d+1}v \in E(T)$. Then $\deg(x_k)=2$ for all $k \in [d+1] \setminus \{1\}$ and $\deg(x_1) = 1$. Consider $T'=T - V(P')$. Then $n'=|V(T')|=n-(d+1)\geq 2(d+1)-d-1=d+1$. Let $D'$ be a $\gamma_d^1(T')$-set. Then $D' \cup \{x_1\}$ is a $d$-distance independent dominating set of $T$. By Corollary~\ref{cor:bipartite-graphs-upper-bound-gamma_d^1},
$$\gamma_d^1(T) \leq |D'|+1 =\gamma_d^1(T')+1 \leq \frac{n'}{d+1} + 1=\frac{n-(d+1)}{d+1}+1=\frac{n}{d+1},$$
and the equality holds if and only if $\gamma_d^1(T) = \gamma_d^1(T')+1$ and $\gamma_d^1(T') =\frac{n'}{d+1}$. The induction hypothesis therefore implies $n'=d+1$ or $T' \in \mathcal{T}_d$.

Suppose $n'=d+1$. Then $n=2(d+1)$ and $T$ is the tree obtained from a copy of $P_{d+1}$ and a tree $T'$ of order $d+1$ by joining $x_{d+1}$ to a vertex $v$ of $T'$. Note that $\diam(T') \leq d$ with equality if and only if $T' \cong P_{d+1}$. Unless $T' \cong P_{d+1}$ and $v$ is a leaf of $T'$, $\{x_{d+1}\}$ is a $d$-distance independent dominating set of $T$, implying that $\gamma_d^1(T)=1<\frac{2(d+1)}{d+1}=\frac{n}{d+1}$, a contradiction. For the exception, we observe $T \cong P_{2(d+1)} \in \mathcal{T}_d$.

Suppose $T' \in \mathcal{T}_d$. Let $T'=T_*' \circ P_d$ for some non-trivial tree $T_*'$. If $v \in V(T_*')$, then $T=T^* \circ P_d \in \mathcal{T}_d$, where $T^*$ is the tree obtained from $T_*'$ by adding a new vertex $x_{d+1}$ and the edge $x_{d+1} v $ to it. If $v \notin V(T_*')$, then let $u_1$ be the $P_d$-support vertex of $T_*'$ such that the attached copy of $P_{d}$ contains $v$. Since $|V(T_*')| \geq 2$, there exists a neighbor $u_2 \in V(T_*')$ of $u_1$. Let $u_1'$ and $u_2'$ be the leaves of $T'$ corresponding to $u_1$ and $u_2$, respectively. Note that $v=u_1'$ is possible, and $D=(L(T') \setminus \{u_1',u_2'\}) \cup \{x_{d+1},u_2\}$ is a $d$-distance independent dominating set of $T$. Thus, $$\gamma_d^1(T) \leq |D|=|L(T')|=|V(T_*')|=\frac{n'}{d+1}
<\frac{n}{d+1}$$
that contradicts our assumption on $T$ and finishes the proof of Claim~\ref{claim:P_d-support-vertex}.
\smallqed

\medskip
Claim~\ref{claim:P_d-support-vertex} shows that if $\gamma_d^1(T)=\frac{n}{d+1}$ and $T$ contains a pendant path $P_{d+1}$, then $T \in \mathcal{T}_d$. The remaining part of the proof verifies that there is no tree $T$ with $|V(T)| > d+1$ and $\gamma_d^1(T)=\frac{n}{d+1}$ that does not contain a pendant $P_{d+1}$.
From now on, we suppose that there is no $P_{d+1}$-support vertex in $T$ and that $\gamma_d^1(T)=\frac{n}{d+1}$. 

Let $s=\diam(T) \geq 2d+1$ and  $P:=v_1v_2 \dots v_{s+1}$ be a diametrical path in $T$. Then $\deg(v_1)=\deg(v_{s+1})=1$. Root $T$ at $v_{s+1}$. Our assumption on the non-existence of  $P_{d+1}$-support vertices and  Claim~\ref{claim:path-support-vertex} imply that the properties stated in Lemma~\ref{lemma:structural} (i)--(v) are valid for $T$.
\medskip
    
    If $s=\diam(T)=2d+1$ then, by Lemma~\ref{lemma:structural} (v), the tree $T$ can be obtained from the $(d-1)$-subdivisions of two stars $K_{1,t_1}$ and $K_{1,t_2}$ with $t_1 \ge t_2 \geq 2$ by joining the centers with an edge.
     Then $N(v_{d+2})$ is a $d$-distance independent dominating set of $T$. Since $d \geq 2$, it gives the following contradiction:
    $$\gamma_d^1(T) \leq |N(v_{d+2})|=t_2+1=\frac{(t_2+1)d+t_2+1}{d+1} < \frac{2dt_2+2}{d+1}\leq \frac{d(t_1+t_2)+2}{d+1}=\frac{n}{d+1}.$$

    So, we may assume that $\diam(T) \geq 2d+2$ and $n \geq 2d+3$. Regarding $v_{d+2}$, we divide the rest of the proof into two cases and prove that in both we get a contradiction.

    \begin{description}
        \item [Case 1.] Each vertex $v$ in $N(v_{d+2}) \setminus \{v_{d+1},v_{d+3}\}$ is of degree at least 3.\\
            By Lemma~\ref{lemma:structural} (iii) and since $P$ is a diametrical path, for each $v \in N(v_{d+2}) \setminus \{v_{d+3}\}$, the subtree $T_v$ is isomorphic to the $(d-1)$-subdivision of a star $K_{1,t_v}$ with $t_v \geq 2$. Clearly, $T_{v_{d+1}}$ is contained in $T_{v_{d+2}}$, and therefore, $|V(T_{v_{d+2}})| \geq 2d+2$. Let $T'=T-T_{v_{d+2}}$. Since $\{v_{d+3}, \dots, v_{2d+3}\} \subseteq V(T')$, we obtain
            $$d+1 \leq n'=|V(T')|  \leq n-2d-2.$$ 

            Let $D'$ be a $\gamma_d^1(T')$-set. Then $D=D' \cup (N(v_{d+2}) \setminus \{v_{d+3}\})$ is a $d$-distance independent dominating set of $T$. Let $p=\deg(v_{d+2})$. Observe that $p \geq 2$ and $n'\leq n-(2d+1) (p-1)-1$. By Corollary~\ref{cor:bipartite-graphs-upper-bound-gamma_d^1}, we get the following contradiction:
            \begin{align*}
		        \gamma_d^1(T) & \leq |D| = \gamma_d^1(T')+p-1\\
                & \leq \frac{n'}{d+1}+p-1\\ 
                &\leq \frac{n-(2d+1)(p-1)-1}{d+1}+p-1\\ 
                &=\frac{n-d(p-1)-1}{d+1}\\
                &<\frac{n}{d+1}.
	        \end{align*}
            
        \item [Case 2.] There is a vertex $v$ in $N(v_{d+2}) \setminus \{v_{d+1},v_{d+3}\}$ with $\deg(v) \leq 2$.\\
        If $\deg(v)=2$ and $T_v$ contains a vertex $u$ with $\deg(u)\ge 3$, then Lemma~\ref{lemma:structural} (iii) implies the existence of a leaf $w \in V(T_u)$ with $d(w,u)=d$. It follows then that $d(w,v_{d+2})\ge d+2$ and $d(w,v_{s+1}) \ge s+1= \diam(T) +1$, a contradiction. Therefore, $\deg(v)\leq 2$ implies that $T_v$ is a path and $v_{d+2}$ is a $P_i$-support vertex for some $i \ge 1$. By our assumption, $i \leq d$.
        Further, by Claim~\ref{claim:path-support-vertex}, we have the following properties.
            \begin{itemize}
                \item If $v_{d+2}$ is a $P_i$-support vertex of $T$ for some $i \in [d-1]$, then there is only one pendant path attached to $v_{d+2}$, and it is clearly of order $i$. 
                \item If $v_{d+2}$ is a $P_d$-support vertex of $T$, then $v_{d+2}$ is not a $P_i$-support vertex of $T$ for any $i \in [d-1]$, and there is at least one copy of $P_d$ attached to $v_{d+2}$.
            \end{itemize}

        \item [Case 2.1.] $L(v_{d+2}) \neq \emptyset$.\\
            In this case $v_{d+2}$ is a $P_1$-support vertex of $T$. Let $x \in L(v_{d+2})$ and $T'=T-x$. Now for each vertex $v \in N(v_{d+2}) \setminus \{v_{d+3}\}$, the subtree $T'_{v}$ is isomorphic to the $P_d$-subdivision of a star $K_{1,t_v}$ for $t_v \geq 2$. Clearly, $n'=|V(T')|=n-1 \geq 2d+2$.
            
            Let $D'$ be a $\gamma_d^1(T')$-set. If $v_{d+2} \in D'$, then $D'$ is also a $d$-distance independent dominating set of $T$. If $v_{d+2} \notin D'$, then since 
            $|D' \cap V(T'_{v_{d+1}})| \geq 1$, we may assume that $v_{d+1} \in D'$. The set $D'$ is also a $d$-distance independent dominating set of $T$. For any subcase, $\gamma_d^1(T) \leq |D'|=\gamma_d^1(T') \leq \frac{n'}{d+1}=\frac{n-1}{d+1}<\frac{n}{d+1}$ by Corollary~\ref{cor:bipartite-graphs-upper-bound-gamma_d^1}.

        \item [Case 2.2.] $L(v_{d+2}) = \emptyset$.\\
            In this case, $v_{d+2}$ is a $P_i$-support vertex of $T$ for some $i \in [d] \setminus \{1\}$ (where if $i = d$, then there could be multiple copies of $P_d$ attached to $v_{d+2}$). Let $P':=x_1x_2\dots x_i$ be the (selected) copy of $P_i$ attached to $v_{d+2}$, where $x_iv_{d+2} \in E(T)$. Then $\deg(x_k)=2$ for all $k \in [i] \setminus \{1\}$ and $\deg(x_1) = 1$. Consider $T'=T-T_{v_d}-T_{x_{i}}$. Then $n'=|V(T')|=n-d-i \leq n-d-2$ and $n' \geq d+3$ since $v_{d+1},v_{d+2}, \dots, v_{2d+3} \in V(T')$.

            Let $D'$ be a $\gamma_d^1(T')$-set. Then $|D' \cap \{v_{d+1},v_{d+2}\}| \leq 1$. If $v_{d+1} \in D'$ and $v_{d+2} \notin D'$, then let $D=D' \cup \{x_i\}$. If $v_{d+1} \notin D'$ and $v_{d+2} \in D'$, then let $D=D' \cup \{v_1\}$. If $v_{d+1},v_{d+2} \notin D'$, then since $v_{d+1}$ is attached to at least two copies of $P_d$, we have $D' \cap (V(T_{v_{d+1}}) \setminus V(T_{v_d})) \neq \emptyset$. Let $D=D' \cup \{v_{d+1},x_i\} \setminus (V(T_{v_{d+1}}) \setminus V(T_{v_d}))$. For any subcase, $D$ is a $d$-distance independent dominating set of $T$, and $\gamma_d^1(T)  \leq |D'|+1=\gamma_d^1(T')+1 \leq \frac{n'}{d+1}+1 \leq \frac{n-d-2}{d+1}+1<\frac{n}{d+1}$ by Corollary~\ref{cor:bipartite-graphs-upper-bound-gamma_d^1}.
   \end{description}
This completes the proof of Theorem~\ref{thm:extremal-trees-gamma_d^1}.
\end{proof}

\section{Upper bounds on $\gamma_d$ and $\gamma_d^1$ of trees in terms of the order and the number of leaves}
\label{sec:upper-bounds-leaves}

For any tree $T$ of order $n$ and with $\ell$ leaves, the set of non-leaves is a dominating set of $T$. Hence, $\gamma_1(T) \leq n-\ell$. Note that the equality holds if and only if each vertex of $T$ is either a leaf or a support vertex. If there exists a vertex $u \in V(T)$ that is neither a leaf nor a support vertex, then $V(T) \setminus (\{u\} \cup L(T))$ is a dominating set of $T$, implying that $\gamma_1(T) < n-\ell$. On the other hand, the upper bound $\gamma_{1}^1(T) \leq n-\ell$ is not true for every tree $T$. For example, let $T' = T^* \circ P_1 \in \mathcal{T}_1$ for some tree $T^*$, and let $T$ be the tree obtained from $T'$ by adding $r \geq 2$ leaves to each vertex of $T'$. It can be checked that $$\gamma_{1}^1(T)=|V(T^*)|+r|V(T^*)|>2|V(T^*)|=2(r+1)|V(T^*)|-2r|V(T^*)|=n-\ell.$$

Set now 
$$\mathcal{F}_{2}  = \left\{T:\ T-L(T)\in \zeta_1 \right\}\,,$$
and if $d\ge 3$, then set 
$$\mathcal{F}_{d}  = \left\{T:\ T-L(T) \text{~is a tree of order~} d \text{~or belongs to~} \mathcal{T}_{d-1}\right\}\,.$$
Note that each graph from $\mathcal{F}_{d}$, $d\ge 2$, is a tree, and the following property is equivalent to the definition of $\mathcal{F}_d$.
\begin{itemize}
    \item[$(\star)$] If $d \ge 3$, a tree $T$ belongs to $\mathcal{F}_{d}$ if and only if it can be obtained from some tree $T'$, which satisfies $|V(T')|=d$ or $T'\in \mathcal{T}_{d-1}$, by adding at least one pendant vertex to each leaf of $T'$, and some number (possibly zero) to other vertices of $T'$. For $d=2$, a tree $T$ belongs to $\mathcal{F}_{2}$ if and only if it can be obtained similarly from a tree $T' \in \zeta_1$. 
\end{itemize}
For $d \geq 2$, we prove the following result.

\begin{theorem}
    \label{thm:trees-upper-bound-leaves-gamma_d^1}
    Let $d \geq 2$ be an integer and $T$ be a tree of order $n$ and with $\ell$ leaves. If $n-\ell \geq d$, then $\gamma_d^1(T) \leq \frac{n-\ell}{d}$ with equality if and only if $T \in \mathcal{F}_d$.
\end{theorem}

\begin{proof}
    Consider the tree $T'=T-L(T)$. Let $n'=|V(T')|=n-\ell \geq d$. Let $D'$ be a $\gamma_{d-1}^1(T')$-set. By Corollary~\ref{cor:bipartite-graphs-upper-bound-gamma_d^1}, $|D'|=\gamma_{d-1}^1(T') \leq \frac{n'}{d}$. Moreover, $D'$ is also a $d$-distance independent dominating set of $T$, implying that 
    \begin{equation} \label{eq:3.7-1}
        \gamma_d^1(T) \leq |D'|=\gamma_{d-1}^1(T') \leq \frac{n'}{d} =\frac{n-\ell}{d}.
    \end{equation}
    
    Assume that $\gamma_d^1(T)  =\frac{n-\ell}{d}$ holds for a tree $T$. Inequalities in (\ref{eq:3.7-1}) therefore imply $\gamma_d^1(T)=\gamma_{d-1}^1(T')=\frac{n'}{d}$. By Theorems \ref{thm:extremal-trees-gamma_1^1} and \ref{thm:extremal-trees-gamma_d^1}, we know that $T' \in \zeta_1$ when $d=2$, and $T'$ is a tree of order $d$ or $T' \in \mathcal{T}_{d-1}$ when $d \geq 3$. Since $T'=T-L(T)$, we conclude $T \in \mathcal{F}_{d}$.

    It remains to prove that $\gamma_d^1(T) \geq \frac{n-\ell}{d}$ holds for every $T \in \mathcal{F}_d$. Consider first a tree $T$ from  $\mathcal{F}_2$ and let $T'= T-L(T)$. Hence $T' \in \zeta_1$. We will prove the inequality by induction on $T'$ according to the recursive definition of $\zeta_1$.
    If $T'\cong K_2$, then $T$ is a double star and $\gamma_2^1(T)= 1 = \frac{n-\ell}{2}$. If $T' \cong D_{r,r}$, for $r \geq 1$, then any $\gamma_1^1(T')$-set is a smallest $2$-distance independent dominating set of $T$, implying that $$\gamma_2^1(T) =\gamma_{1}^1(T')=r+1 =\frac{2r+2}{2}= \frac{n'}{2} =\frac{n-\ell}{2}.$$
    
    Assume next that $T'=T-L(T)$ is a tree from $\zeta_1$ which is neither $K_2$ nor a double star. Let $T_2'=T'$ and let $T_1'$ be the tree from $\zeta_1$ such that $T_2'$ is obtained from $T_1'$ by the recursive construction of $\zeta_1$, that is, $T_2'$ can be obtained by joining the center $u$ of a new copy of $K_{1,t}\ (t\geq 1)$ to a support vertex $v$ of $T_1'$, and adding $t-1$ leaves at $v$. For $i \in [2]$, let $T_i$ be a tree from $\mathcal{F}_2$, which is obtained from $T_i'$ according to $(\star)$. Moreover, let $n_i'=|V(T_i')|$, $n_i=|V(T_i)|$, and $\ell_i=|L(T_i)|$, $i \in [2]$. 
    
    Assume that $\gamma_2^1(T_1) = \frac{n_1-\ell_1}{2}$. We are going to  
    prove that $\gamma_2^1(T_2) \geq \frac{n_2-\ell_2}{2}$. Note that $n_i'=n_i-\ell_i$ and $n_2'=n_1'+2t$. Let $D_2$ be a $\gamma_{2}^1(T_2)$-set that contains as few leaves from $T_2$ as possible and let $D_1=D_2 \cap V(T_1)$. If $v \in D_2$, then $u \notin D_2$ and, by the minimality of $|D_2 \cap L(T_2)|$, we have $L_{T_2'}(u) \subset D_2$. Now $D_1=D_2 \setminus L_{T_2'}(u)$ is a $2$-distance independent dominating set of $T_1$, implying that $\gamma_2^1(T_1) \leq |D_1| =|D_2|-t$. If $v \notin D_2$, then we may assume that $u \in D_2$. Also, $L_{T_2'}(v) \setminus L_{T_1'}(v) \subset D_2$ holds by the minimality of $|D_2 \cap L(T_2)|$. Further, $\emptyset \neq L_{T_1'}(v) \subset D_1$, and $v$ and the leaves added to $L_{T_1'}(v)$ in $T_1$ will be independently dominated by $L_{T_1'}(v)$. Hence, $D_1=D_2\setminus \{u\} \setminus (L_{T_2'}(v) \setminus L_{T_1'}(v))$ is a $2$-distance independent dominating set of $T_1$, implying that $\gamma_2^1(T_1) \leq |D_1|=|D_2|-t$. Hence no matter whether $v$ belongs to $D_2$ or not, we have
    $$\gamma_2^1(T_2) \geq \gamma_2^1(T_1)+t = \frac{n_1-\ell_1}{2}+\frac{n_2'-n_1'}{2}=\frac{n_2'}{2}=\frac{n_2-\ell_2}{2}.$$
    
    Assume now that $T \in \mathcal{F}_d$ and $d \ge 3$. For $T'=T-L(T)$, let $n'=|V(T')|=n-\ell \geq d$. If $n'=d$, then $\gamma_d^1(T) \geq 1=\frac{n-\ell}{d}$. If $T' \in \mathcal{T}_{d-1}$, then let $T'=T^* \circ P_{d-1}$ for some non-trivial tree $T^*$. For each $u \in V(T^*)$, let $T_u'$ be the subtree of $T'$ induced by $u$ and the vertices of the copy of $P_{d-1}$ attached to $u$, and let $T_{u}$ be the subtree of $T$ induced by $V(T_u')$ and the leaves added to $V(T_u')$ in $T$.   
    If $D$ is a $\gamma_d^1(T)$-set, then $|D \cap V(T_u)| \geq 1$ for every $u \in V(T^*)$. Thus, we have $$\gamma_d^1(T) =|D| \geq |V(T^*)|=\frac{n'}{d} = \frac{n-\ell}{d}.$$
    This completes the proof of Theorem~\ref{thm:trees-upper-bound-leaves-gamma_d^1}.
\end{proof}

We note that the condition of $n \geq d+\ell$ is necessary in Theorem~\ref{thm:trees-upper-bound-leaves-gamma_d^1}. Let $T'$ be a tree of order at most $d-1$. Consider the tree $T$ obtained from $T'$ by adding at least one pendant vertex to each leaf of $T'$ and some number to other vertices of $T'$. Then $n'=|V(T')|=n-\ell \leq d-1$ and we may infer $\gamma_d^1(T) \geq 1 > \frac{d-1}{d} \geq \frac{n-\ell}{d}$.
\medskip

Favaron~\cite{Favaron-1992} proved that if $T$ is a tree of order $n \geq 2$ and with $\ell$ leaves, then $\gamma_1^1(T) \leq \frac{n+\ell}{3}$, and gave the full list of extremal trees for this bound. Our next theorem extends Favaron's result to all $d \geq 2$.

\begin{theorem}
    \label{thm:trees-upper-bound+leaves-gamma_d^1}
    Let $d \geq 2$ be an integer and $T$ a tree of order $n$ and with $\ell$ leaves. If $n \geq d$, then $\gamma_d^1(T) \leq \frac{n+\ell}{d+2}$ with equality if and only if $T \in \{P_d\} \cup \mathcal{T}_d$. 
\end{theorem}

\begin{proof}
    If $T \cong P_d$, then $\gamma_d^1(T)=1=\frac{d+2}{d+2}=\frac{n+\ell}{d+2}$. If $T \in \mathcal{T}_d$, then by Proposition \ref{prop:extremal-bipartite-graphs-value-gamma_d^1}, $\gamma_d^1(T)=\frac{n}{d+1}=\frac{n+\frac{n}{d+1}}{d+2}=\frac{n+\ell}{d+2}$. To prove the upper bound and that the equality implies $T \in \{P_d\} \cup \mathcal{T}_d$, we proceed by induction on $n$. If $\diam(T) \leq 2d$, then $\gamma_d^1(T) =1 = \frac{d+2}{d+2} \leq \frac{n+\ell}{d+2}$. The equality holds if and only if $n=d$ and $\ell=2$, implying that $T \cong P_d$. So, we may assume that $\diam(T) \geq 2d+1$ and $n \geq 2d+2$. Note that if $\ell >\frac{n}{d+1}$, then by Corollary~\ref{cor:bipartite-graphs-upper-bound-gamma_d^1}, $\gamma_d^1(T) \leq \frac{n}{d+1}< \frac{n+\ell}{d+2}$. 

    \begin{claim}
        \label{claim:path-support-vertex-n+l}
        Let $i \in [d-1]$ and $j \in [d]$ with $i \leq j$. If $T$ has a vertex $v$ that is a $(P_i,P_j)$-support vertex, then $\gamma_d^1(T) < \frac{n+\ell}{d+2}$.
        
    \end{claim}

    \proof
        Let $P':=x_1x_2\dots x_i$ and $P'':=y_1y_2\dots y_j$ be a copy of $P_i$ and $P_j$, respectively, attached to $v$ in $T$, where $x_iv,y_jv \in E(T)$. Since $n \geq 2d+2$ and $|V(P')| \leq |V(P'')| \leq d$, we have $\deg(v) \geq 3$. 
        Consider $T'=T-V(P')$. Then $\ell'=|L(T')| = \ell-1$ and 
        $n'=|V(T')|=n-i \geq d+3$. As in the proof of Claim~\ref{claim:path-support-vertex} it can be proved that there exists a $\gamma_d^1(T')$-set $D'$ that is a $d$-distance independent dominating set of $T$. Using the induction hypothesis, we have 
        $$\gamma_d^1(T) \leq |D'|=\gamma_d^1(T') \leq \frac{n'+\ell'}{d+2}= \frac{n-i+\ell-1}{d+2} 
        <\frac{n+\ell}{d+2}.$$
        This proves Claim \ref{claim:path-support-vertex-n+l}.
    \smallqed
    
    \begin{claim}
        \label{claim:P_d-support-vertex-n+l}
        If $T$ has a $P_{d+1}$-support vertex $v$, then $\gamma_d^1(T) \leq \frac{n+\ell}{d+2}$ and if equality holds, then $T \in \mathcal{T}_d$.
    \end{claim}

    \proof
        Let $P':=x_1x_2\dots x_{d+1}$ be a copy of $P_{d+1}$ attached to $v$, where $x_{d+1}v \in E(T)$. Then $\deg_T(x_k)=2$ for all $k \in [d+1] \setminus \{1\}$ and $\deg_T(x_1) = 1$. Consider $T'=T - V(P')$. Since $n \geq 2d+2$, we have $\deg_T(v) \geq 2$. Then $\ell'=|L(T')| =\ell$ if $\deg_T(v)=2$ and $\ell'=\ell-1$ if $\deg_T(v)\geq 3$. 
        We observe that $n'=|V(T')|=n-(d+1)\geq d+1$ and consider two cases according to the degree of $v$.
        \begin{description}
        \item [Case D1.] $\deg_T(v)\geq 3$.\\
        Let $D'$ be a $\gamma_d^1(T')$-set. The set $D' \cup \{x_1\}$ is a $d$-distance independent dominating set of $T$. By the induction hypothesis, $$\gamma_d^1(T) \leq |D' \cup \{x_1\}| =\gamma_d^1(T')+1 \leq \frac{n'+\ell'}{d+2} + 1=\frac{n-(d+1)+\ell-1}{d+2}+1=\frac{n+\ell}{d+2},$$
        and the equality holds if and only if $\gamma_d^1(T) = \gamma_d^1(T')+1$ and $\gamma_d^1(T') =\frac{n'+\ell'}{d+2}$. Note that $n' \geq d+1$, so $T' \ncong P_d$ and $T' \in \mathcal{T}_d$.

        Let $T'=T_*' \circ P_d$ for some non-trivial tree $T_*'$. Then $\ell'=\frac{n'}{d+1}$. Since $\deg(v) \geq 3$, we infer that $v \notin L(T')$. If $v \in V(T_*')$, then $T=T_* \circ P_d \in \mathcal{T}_d$, where $T_*$ is the tree obtained from $T_*'$ by adding a new vertex $x_{d+1}$ to it such that $x_{d+1} v \in E(T_*)$. If $v \notin V(T_*') \cup L(T')$, then let $u$ be the $P_d$-support vertex of $T_*'$ attached to the copy of $P_{d}$ containing $v$, and $u'$ be the leaf of $T'$ corresponding to $u$. Note that $v \neq u'$ and $D=(L(T') \setminus \{u'\}) \cup \{x_{d+1}\}$ is a $d$-distance independent dominating set of $T$. Thus, $$\gamma_d^1(T) \leq |D|=|L(T')|=\frac{n'}{d+1}=\frac{n'+\ell'}{d+2}=\frac{n-(d+1)+\ell-1}{d+2}<\frac{n+\ell}{d+2}.$$
        \item [Case D2.] $\deg_T(v)=2$.\\
          Let $P'':=x_1x_2\dots x_{d+1}v$ be a copy of $P_{d+2}$ attached to $v'$, where $vv' \in E(T)$. Consider $T''=T-V(P'')=T'-v$. Then $n''=|V(T'')|=n-(d+2) \geq d$, and $\ell''=|L(T'')| \leq \ell$ with equality if and only if $\deg_T(v')=2$. Let $D''$ be a $\gamma_d^1(T'')$-set. The set $D'' \cup \{x_2\}$ is a $d$-distance independent dominating set of $T$. By the induction hypothesis,
        $$\gamma_d^1(T) \leq |D'' \cup \{x_2\}| =|D''|+1 =\gamma_d^1(T'')+1 \leq \frac{n''+\ell''}{d+2} + 1 \leq \frac{n-(d+2)+\ell}{d+2}+1=\frac{n+\ell}{d+2},$$
        and the equality holds if and only if $\gamma_d^1(T) = \gamma_d^1(T'')+1$, $\ell''=\ell$, and $\gamma_d^1(T'') =\frac{n''+\ell''}{d+1}$, i.e., $T'' \in \{P_d\} \cup \mathcal{T}_d$.

        Note that $\deg_T(v')=2$ and $\deg_{T''}(v')=1$. If $T'' \cong P_d$, then $T \cong P_{2d+2} \in \mathcal{T}_d$. Suppose that $T'' \in \mathcal{T}_d$. Let $T''=T_*'' \circ P_d$ for some non-trivial tree $T_*''$. Then $\ell''=\frac{n''}{d+1}$. Clearly, $v' \in L(T'')$. Let $u_1' \in V(T_*'')$ be the $P_d$-support vertex in $T''$, which is attached to the copy of $P_{d}$ containing $v'$. Since $|V(T_*'')| \geq 2$, there exists a neighbor $u_2' \in V(T_*'')$ of $u_1'$. It is clear that $v'$ is the leaf of $T'$ corresponding to $u_1'$. Let $u_2''$ be the leaf of $T'$ corresponding to $u_2'$. Since $d \geq 2$, the set $D=(L(T'') \setminus \{v',u_2''\}) \cup \{u_2',x_{d+1}\}$ is a $d$-distance independent dominating set of $T$. Thus, we have
        $$\gamma_d^1(T) \leq |D|=|L(T'')|=\frac{n''}{d+1}=\frac{n''+\ell''}{d+2}=\frac{n-(d+2)+\ell}{d+2}<\frac{n+\ell}{d+2}.$$
        \end{description}
        This completes the proof of Claim \ref{claim:P_d-support-vertex-n+l}.
    \smallqed

 \medskip
   In the continuation, we may suppose that there is no $P_{d+1}$-support vertex in $T$ and also that if $v$ is a $(P_i,P_j)$-support vertex, then $i=j=d$. Let $s=\diam(T) \geq 2d+1$ and let $P:=v_1v_2 \dots v_{s+1}$ be a diametrical path in $T$. Root $T$ at $v_{s+1}$. Hence, by Lemma~\ref{lemma:structural},
    $\deg(v_k)\leq 2$ for each $k \in [d] \cup ([s+1] \setminus [s-d+1])$, and   $\deg(v_{k}) \geq 3$ for each $k \in \{ d+1, s - d + 1 \}$. It also follows that the subtree $T_{v_{d+1}}$ is isomorphic to the $(d-1)$-subdivision of a star $K_{1,t}$ for some $t \geq 2$.
    \medskip

    If $s=\diam(T)=2d+1$, then by Lemma~\ref{lemma:structural} (v), $T$ is obtained from the $(d-1)$-subdivision of a star 
    $K_{1,t_1}$ and the $(d-1)$-subdivision of a star
    $K_{1,t_2}$ by joining the centers $v_{d+1}$ and $v_{d+2}$. We may assume that $t_1 \geq t_2 \geq 2$. Then $N(v_{d+2})$ is a $d$-distance independent dominating set of $T$. Since $d \geq 2$, we have 
    \begin{align*}
            \gamma_d^1(T) & \leq |N(v_{d+2})|=\deg(v_{d+2})=t_2+1 = \frac{(d+1)t_2+d+t_2+2}{d+2} \\ & < \frac{(d+1)t_2+dt_2+t_2+2}{d+2} = \frac{2(d+1)t_2+2}{d+2} \\
            & \leq \frac{d(t_1+t_2)+2+(t_1+t_2)}{d+2}=\frac{n+\ell}{d+2}.
        \end{align*}
        
        So, we may assume that $\diam(T) \geq 2d+2$ and $n \geq 2d+3$.  Regarding $v_{d+2}$, we divide the rest of the proof into two cases and prove that the strict inequality $\gamma_d^1(T) < \frac{n+\ell}{d+2}$ holds in each case.

    \begin{description}
        \item [Case 1.] Every vertex $v$ in $N(v_{d+2}) \setminus \{v_{d+1},v_{d+3}\}$ is of degree at least 3.\\
            For each vertex $v \in N(v_{d+2}) \setminus \{v_{d+3}\}$ we have $\deg(v) \geq 3$, and the subtree $T_v$ is isomorphic to the $(d-1)$-subdivision of a star $K_{1,t_v}$ for $t_v \geq 2$. Let $T'=T-T_{v_{d+2}}$ and $p= \deg(v_{d+2})$. 
            It holds that
            $$d+1 \leq n'=|V(T')| \leq n-1- (2d+1) (p-1).$$
             Moreover, we have $$\ell'=|L(T')| 
             \leq \ell-2(p-1)+1=\ell-2p-1,$$ with equality if and only if $\deg(v_{d+3})=2$, and for each $v \in N(v_{d+2}) \setminus \{v_{d+3}\}$, $t_v=2$.
            
            Let $D'$ be a $\gamma_d^1(T')$-set. Then $D=D' \cup (N(v_{d+2}) \setminus \{v_{d+3}\})$ is a $d$-distance independent dominating set of $T$.  Since $d \geq 2$ and $p \geq 2$, by the induction hypothesis we get
            \begin{equation*}
	        \begin{aligned}
		        \gamma_d^1(T) \leq |D|&=|D'|+p-1=\gamma_d^1(T')+p-1 \leq \frac{n'+\ell'}{d+2}+p-1\\ &\leq \frac{n-1-(2d+1)(p-1)+\ell-2p-1}{d+2}+p-1\\ &=\frac{n+\ell-d(p-1)-p-3}{d+2}<\frac{n+\ell}{d+2}.
	        \end{aligned}
            \end{equation*}

        \item [Case 2.] There is a vertex $v$ in $N(v_{d+2}) \setminus \{v_{d+1},v_{d+3}\}$ with $\deg(v) \leq 2$.\\
            Since $v_{d+2}$ is not a $P_{d+1}$-support vertex and $P$ is a diametrical path, Lemma~\ref{lemma:structural} (iii) implies that $T_v$ is a pendant path $P_i$ for some $i \in [d]$. Moreover, we have the following.
            \begin{itemize}
                \item If $v_{d+2}$ is a $P_i$-support vertex of $T$ for some $i \in [d-1]$, then there is no other pendant path attached to $v_{d+2}$.
                
                \item If $v_{d+2}$ is a $P_d$-support vertex of $T$, then $v_{d+2}$ is not a $P_i$-support vertex of $T$ for any $i \in [d-1]$, and there is at least one copy of $P_d$ attached to $v_{d+2}$. 
            \end{itemize}

        \item [Case 2.1.] $L(v_{d+2}) \neq \emptyset$.\\
            Let $x \in L(v_{d+2})$ and $T'=T-x$. Clearly, $\deg_T(v_{d+2}) \geq 3$ and $\deg_{T'}(v_{d+2})\geq 2$. Then $\ell'=|L(T')|=\ell-1$ 
            and
            $n' =|V(T')|=n-1 \geq 2d+2$. Let $D'$ be a $\gamma_d^1(T')$-set. By considering whether $v_{d+2}$ is in $D'$ or not, we observe that $D'$ can be chosen such that it is also a  $d$-distance (independent) dominating set of $T$. By the induction hypothesis, we have $\gamma_d^1(T) \leq |D'|=\gamma_d^1(T') \leq \frac{n'+\ell'}{d+2}
            =\frac{n-1+\ell-1}{d+2}
            <\frac{n+\ell}{d+2}$.

        \item [Case 2.2.] $L(v_{d+2}) = \emptyset$.\\
            Let $P':=x_1x_2\dots x_i$ be a copy of $P_i$ attached to $v_{d+2}$, where $x_iv_{d+2} \in E(T)$. Then $i \in [d]\setminus \{1\}$, and $\deg(x_k)=2$ for all $k \in [i] \setminus \{1\}$, while $\deg(x_1) = 1$. Consider $T'=T-T_{v_d}-T_{x_{i}}$. By Lemma~\ref{lemma:structural} (ii), $\deg(v_{d+1}) \geq 3$ and, by our condition, $\deg(v_{d+2}) \geq 3$. Therefore, $\ell'=|L(T')|=\ell-2$. We also know that $n'=|V(T')|=n-d-i \leq n-d-2$ and $n' \geq d+3$.

            Let $D'$ be a $\gamma_d^1(T')$-set. Then $|D' \cap \{v_{d+1},v_{d+2}\}| \leq 1$. As in Case 2.2 of Theorem \ref{thm:extremal-trees-gamma_d^1}, let
            $$
            D = 
                \begin{cases}
	            D' \cup \{x_i\}, &\text{if~} v_{d+1} \in D' \text{~and~} 
                        v_{d+2} \notin D' ,\\
	              D' \cup \{v_1\}, &\text{if~} v_{d+1} \notin D' 
                        \text{~and~} v_{d+2} \in D' ,\\
	              D' \cup \{v_{d+1},x_i\} \setminus (V(T_{v_{d+1}}) 
                        \setminus V(T_{v_d})), &\text{if~} v_{d+1},v_{d+2} \notin D' .\\
                \end{cases}
            $$
            For any subcase, $D$ is a $d$-distance independent dominating set of $T$. By the induction hypothesis, we have $\gamma_d^1(T) \leq |D| \leq |D'|+1=\gamma_d^1(T')+1 \leq \frac{n'+\ell'}{d+2}+1 \leq \frac{n-d-2+\ell-2}{d+2}+1<\frac{n+\ell}{d+2}$.
    \end{description}
This completes the proof of Theorem \ref{thm:trees-upper-bound+leaves-gamma_d^1}.
\end{proof}

Now we set
$$\mathcal{F}_{2}'  = \left\{T:\ T-L(T) \in \{K_2\} \cup \mathcal{T}_{1}\right\}\,,$$
and if $d \geq 3$, then set $$\mathcal{F}_{d}'=\mathcal{F}_{d}.$$

By Theorems \ref{thm:trees-upper-bound-leaves-gamma_d^1} and \ref{thm:trees-upper-bound+leaves-gamma_d^1}, we have the following two corollaries, respectively.

\begin{corollary}
    \label{cor:trees-upper-bound-leaves-gamma_d^0}
    Let $d \geq 2$ be an integer and $T$ be a tree of order $n$ and with $\ell$ leaves. If $n-\ell \geq d$, then $\gamma_d(T) \leq \frac{n-\ell}{d}$ with equality if and only if $T \in \mathcal{F}_d'$.
\end{corollary}

\begin{corollary}
    \label{cor:trees-upper-bound+leaves-gamma_d^0}
    Let $d \geq 2$ be an integer and $T$ be a tree of order $n$ and with $\ell$ leaves. If $n \geq d$, then $\gamma_d(T) \leq \frac{n+\ell}{d+2}$ with equality if and only if $T \in \{P_d\} \cup \mathcal{T}_d$.
\end{corollary}

Combining the above results with Corollary~\ref{cor:bipartite-graphs-upper-bound-gamma_d^1}, we obtain

\begin{corollary}
    \label{cor:piecewise-function}
    If $d\ge 2$, and $T$ is a tree with $\ell$ leaves and of order $n \geq d+\ell$, then
    $$
    \gamma_d(T) \leq \gamma_d^1(T) \leq 
    \begin{cases}
	\frac{n-\ell}{d}, &\text{if~} n<(d+1)\ell ,\\
	\frac{n}{d+1}, &\text{if~} n=(d+1)\ell ,\\
	\frac{n+\ell}{d+2}, &\text{if~} n>(d+1)\ell .\\
    \end{cases}
    $$
    Moreover, these bounds are best possible.
\end{corollary}

\section{A conjecture}

Recall that Ma and Chen~\cite{MaDeXing-2004} described equivalently bipartite graphs $G$ of order $n$ with $\gamma_1^1(G)=\frac{n}{2}$. For $d \geq 2$ we pose: 

\begin{conjecture}
\label{conj:bip}
If $d\ge 2$ and $G$ is a connected bipartite graph $G$ of order $n$, then $\gamma_d^1(G)=\frac{n}{d+1}$ if and only if $G \in \{C_{2d+2}\} \cup \mathcal{B}_d$ or $n=d+1$.     
\end{conjecture}

Since $\gamma_1^1(K_{r,r})=r=\frac{r+r}{2}=\frac{n}{2}$, the condition of $d \geq 2$ of the conjecture above is necessary. If Conjecture~\ref{conj:bip} holds true, then it generalizes Theorem~\ref{thm:extremal-trees-gamma_d^1}. Moreover, the result~\cite[Theorem 3]{Topp-1991} due to Topp and Volkmann, restricted to bipartite graphs, gives exactly the same characterization for graphs $G$ with $\gamma_d(G) = \frac{n}{d+1}$ as we pose in Conjecture~\ref{conj:bip} for the $d$-distance independent domination. 

\section*{Acknowledgements}

Csilla Bujt\'{a}s, Vesna Ir\v{s}i\v{c}, and Sandi Klav\v{z}ar were supported by the Slovenian Research and Innovation Agency (ARIS) under the grants P1-0297, N1-0285, N1-0355, and  Z1-50003.  Vesna Ir\v{s}i\v{c} also acknowledges the financial support from the European Union (ERC, KARST, 101071836).

\end{document}